\documentclass{article}

\usepackage{PRIMEarxiv}

\usepackage[utf8]{inputenc} 
\usepackage[T1]{fontenc}    
\usepackage{amsmath,amsfonts} 
\usepackage{algorithmic}
\usepackage{algorithm}  
\usepackage{array} 
\usepackage[caption=false,font=normalsize,labelfont=sf,textfont=sf]{subfig}
\usepackage{textcomp}
\usepackage{stfloats}
\usepackage{url}
\usepackage{verbatim} 
\usepackage{graphicx}
\usepackage{lipsum}
\usepackage{cite}
\usepackage{amsthm}
\usepackage{pgfplots} 
\pgfplotsset{compat=1.6}
\usepackage{amssymb}
\usepackage{csvsimple}
\usepackage{tikz} 
\usepackage{tabularx}
\usepackage{booktabs} 
\usetikzlibrary{shapes,arrows,3d,positioning,calc}  
\newtheorem{remark}{Remark} 
 
\newtheorem{prop}{Proposition}  

\newtheorem{lemma}{Lemma} 

\DeclareMathOperator*{\argmin}{arg\,min}
\newcolumntype{C}{>{\centering\arraybackslash}X} 
\usepackage{xcolor}

\title{$\ell^1$-norm Regularized $\ell^1$-norm Best-fit Lines}

\author{
  Xiao Ling\\  
  Virginia Commonwealth University \\
  \texttt{lingx@vcu.edu} \\ 
   \And
   Paul Brooks\\ 
  Virginia Commonwealth University \\
  Richmond\\
  \texttt{jpbrooks@vcu.edu} \\ 
}

\begin{document} 
\maketitle

\begin{abstract} 
This work develops a sparse and outlier-insensitive method to fit a one-dimensional subspace that can be used as a replacement for eigenvector methods such as principal component analysis (PCA). The method is insensitive to outlier observations by formulating procedures as optimization problems that seek the best-fit line according to the $\ell^1$ norm. It is also capable of producing sparse principal components by leveraging an additional penalty term induce sparsity. The algorithm has a worst-case time complexity of $O{(m^2n \log n)}$ and, under certain conditions, produces a globally optimal solution. An implementation of this algorithm in the parallel and heterogeneous environment NVIDIA CUDA is tested on synthetic and real world datasets and compared to existing methods. The results demonstrate the scalability and efficiency of the proposed approach. 
\end{abstract}

\section{Introduction} 
Subspace estimation can be used for dimension reduction by projecting data in a high-dimensional feature space to a low-dimensional subspace. It sheds light on a broad range of tasks from computer vision to pattern recognition. Conventional PCA, hereafter referred to as $\ell^2$-PCA, is a widely used technique to finding a best-fit subspace. $\ell^2$-PCA produces the linear combinations of the original features such that the combinations capture maximal variance. $\ell^2$-PCA can be computed via the singular value decomposition (SVD) of the data matrix. The $\ell^2$ metric is sensitive to outliers in the data matrix. A solution to this disadvantage is replacing the $\ell^2$ metric with an $\ell^1$-norm analog \cite{Kwak2008PrincipalMaximization,Candes2011RobustAnalysis,Brooks2017IdentifyingRanking,Markopoulos2018AdaptiveRejection,Ke2003RobustNorm}. When the feature space has large number of dimensions, it can be difficult to interpret the principal component (PC) loadings and assess important input features.  

To help with interpretation, we can encourage sparsity in the PC loadings. Many methods have been developed that apply $\ell^1$-regularization since \cite{Tibshirani1996RegressionLasso} proposed the LASSO method for regression problems. \cite{Wang2007RobustLAD-lasso} demonstrated the efficacy of $\ell^1$-regularization as a vehicle for inducing sparsity. A simple and intuitive definition of sparsity of data is the number of nonzero entries in the dataset, quantified by $\ell^0$ norm. 

In this chapter, we propose an algorithm with modest computational requirements for $\ell^1$ regularization with the traditional squared $\ell^2$-norm error replaced by the $\ell^1$ measure. Consider the optimization problem to find an $\ell^1$-norm regularized $\ell^1$-norm best-fit one-dimensional subspace: 
\begin{align}
\label{formulation1}
    &\min_{v,\alpha} \sum_{i\in N}\|x_i - v\alpha_i\|_1 + \lambda\|v\|_1,
\end{align}
where $x_i$, $i \in N$ are points in $\mathbb{R}^m$. An optimal vector $v^*$ determines a line through the origin corresponding to the best-fit subspace. For each point $x_i$, the optimal coefficient $\alpha_i^*$ specifies the locations of the projected points $v^*\alpha_i$ on the line defined by $v^*$. Due to the nature of the $\ell^1$ norm, some components of $v$ will be reduced to zero if $\lambda$ is large enough. Therefore, our proposed method simultaneously generates both a best-fit and a sparse line in $m$ dimensions, which makes it suitable for large or high-dimensional data. The method can be extended to the problem of fitting subspaces. The problem in \eqref{formulation1} is non-linear, non-convex, and non-differentiable. Therefore, we adapt the approximation algorithm of \cite{Brooks2017IdentifyingRanking} to the regularized problem. 
\section{Related Works} 
Subspace fitting using the $\ell^1$ norm has a long history including the work of Boscovitch in the 1760s and Edgeworth and Laplace in the 1880s \cite{bloomfield} on fitting a line in two dimensions.   \cite{Ke2003RobustNorm,Ke2005RobustProgramming} present an alternative minimization algorithm using weighted medians and convex quadratic programming with random initialization. \cite{brooks2013pure} describe a polynomial-time algorithm for finding an $\ell^1$-norm best-fit hyperplane using linear programming (LP). \cite{Song2017LowError} propose a polynomial-time algorithm to approximate an $\ell^1$-norm low-rank subspace. \cite{brooks2017estimating} demonstrate an equivalence between their approach, that of \cite{Tsagkarakis2017OnRank}, and that of \cite{Brooks2017IdentifyingRanking}.  \cite{gillis2018} recently showed that finding an $\ell^1$-norm best-fit line is NP-hard.
The best-fit subspace problem is closely related to the low-rank approximation problem. There are growing needs for a robust and sparse alternative to PCA, and this has led to an active research area in low-rank approximation. In this section, we will describe the general idea of the low-rank approximation and some popular methods for solving the problems. 

The nuclear norm of $L$, $\|L\|_*$, is the sum of the singular values of $L$. Minimizing $\|L\|_*$ encourages $L$ to be low rank. Low-rank approximation is based on the assumption that a matrix $X$ can be decomposed into a low-rank matrix $L$ and a sparse matrix $S$; namely, $X=L+S$. Minimizing the $\ell^1$ norm of $S$ encourages sparsity. For a matrix $X\in\mathcal{R}^{n\times m}$, $\|X\|_{2,1}= \sum_{i=1}^n\sqrt{\sum_{j=1}^m X_{ij}^2}$. The $\|X\cdot\|_{2,1}$ can force $L$ to have zero columns corresponding to outliers. 

Given a collection of centered data points $X\in\mathbb{R}^{n\times m}$, the low-rank approximation problem can be formulated as follows.
\begin{align}
    \min_{L} \|X-L\|_p\qquad \mbox{s.t. }\ rank(L)\leq r  \nonumber
\end{align}
where $\|\cdot\|_p$ can be the $\ell^1$ or $\ell^2$ norm. The problem is nonconvex because of the rank calculation. The low-rank approximation problem Principal Component Pursuit (PCP) problem studied by \cite{chandrasekaran2011rank,Candes2011RobustAnalysis} and Low-Rank and Block-Sparse Matrix Decomposition (LRBS) \cite{tang2011robust} are formulated as:     
\begin{align}
     \textsc{PCP}: & \min_{L,S} & \|L\|_*+\lambda \|S\|_1 \nonumber\\ 
     & \mbox{s.t. } &  L+S=X.   \label{eq:pcp}\\
     \textsc{LRBS}: & \min_{L,S} & \|L\|_*+k\lambda \|S\|_{2,1}+k(1-\lambda)\|L\|_{2,1}, \nonumber \\ 
     & \mbox{s.t. } & L+S=X.   \label{eq:lrbs}
\end{align}
where $L$ is low-rank matrix and $S$ is the sparse matrix. Problems \eqref{eq:pcp} and \eqref{eq:lrbs} are convex and can be solved by various methods. An augmented Lagrange multiplier (ALM) algorithm is used in \cite{Candes2011RobustAnalysis} and an alternating splitting augmented Lagrangian method (ASALM) algorithm is used in \cite{adm}. 
For these and related methods, a rank-one component can be obtained by the SVD of $L$, which will be used to compare with our proposed method in computational experiments.

ALM and ADM methods use the following Lagrangian function of \eqref{eq:pcp}:
\begin{align}
   \mathcal{L}(L,S,Z)=\|L\|_*+\lambda \|S\|_1+Z^T(X-L-S)+\frac{\beta}{2}\|X-L-S\|_F^2,  \label{eq:lag1}  
\end{align}
where $\beta$ is a positive penalty parameter. The difference between ALM and ADM is that ALM minimizes \eqref{eq:lag1} with respect to $L$ and $S$, by setting $(L^{j+1},S^{j+1})=\argmin \mathcal{L}(L,S,Z^j)$ and updating the Lagrangian multiplier matrix $Z^{j+1}=Z^j+\beta(X-L^j-S^j)$\cite{liu2012robust}. \cite{min2010decomposing} introduces a related variant ALM scheme. In contrast, ADM sets $L^{j+1}=\argmin\mathcal{L}(L,S^j,Z^j)$, $S^{j+1}=\argmin\mathcal{L}(L^{j},S,Z^j)$ and $Z^{j+1}=Z^j+\beta(X-L^{j}-S^j)$. The form of $L^{j+1}$ and $S^{j+1}$ usually have closed form solutions \cite{yang2011alternating,liu2010robust,cai2010singular,lin2011linearized}. \eqref{eq:lrbs} can be solved by ALM by minimizing the following Lagrangian function of \eqref{eq:lrbs}\cite{tang2011robust}:
\begin{align}
    \mathcal{L}(L,S,Z)=&\|L\|_*+k\lambda \|S\|_{2,1}+k(1-\lambda)\|L\|_{2,1}+ \nonumber \\
   & Z^T(X-L-S)+\frac{\beta}{2}\|X-L-S\|_F^2.  \label{eq:lag2}   
\end{align} 

The problem \eqref{eq:pcp} is convex and can be solved by ALM-based algorithms. \cite{lee2016computationally} developed an updating scheme with closed-form solutions at each ALM iteration. An ALM variant called DNDP-ALM proposed in \cite{cheng2017double} improved the original optimization problem and incorporated noise into the constraints. \cite{bai2018generalized,shen2020partial} solves PCP by taking advantages of the multi-block structure. \cite{lin2010augmented} introduce the exact and inexact ALM methods that achieves a good performance in solving the PCP problem. \cite{RodriguezPaul2013Fpcp} proposed an algorithm approximately an
order of magnitude faster than inexact ALM to construct a sparse component of the same quality. \cite{ShenYuan2019Aamm} described a simple and almost parameter-free algorithm by reformulating the PCP as an unconstrained nonconvex program and then performing alternating minimization scheme. A parallel splitting ALM method was introduced in \cite{liu2020parallel}. \cite{han2013adm} described an ADM algorithm is able to achieve global convergence under standard assumptions. \cite{liu2011solving} proposes the first linear time algorithm for exactly solving very large PCP problems. A scalable algorithm proposed in \cite{6797367} is able to generate suboptimal solution to PCP. \cite{he2012alternating} combines an ADM with a Gaussian back substitution procedure \cite{zheng2012practical} to solve the PCP. 

The LRSLibrary \cite{sobral2016lrslibrary} provides various implementations of algorithms and variants to the PCP problem in MATLAB.

\cite{de2001robust} presented a method for robust principal component analysis (RPCA) that can be used for automatic learning of subspace for data. \cite{6738015} proposed a simple alternating minimization algorithm for
solving a minor variation on the original Principal Component Pursuit (PCP). Under the same assumption to PCP, similar problem formulations have been studied. \cite{NIPS2010} described a problem $\min_{L,S} \|L\|_*+\lambda \|S\|_{1,2}$ called Outlier Pursuit can be efficiently solved by proximal gradient algorithm \cite{cai2010singular}. \cite{kang2015robust} presented a more robust and less biased nonconvex formulation and solved using augmented Lagrange multiplier framework. \cite{zhou2011godec} described a novel low-rank and sparse decomposition problem called Go Decomposition (GoDec).  

On the other hand, $\ell^1$-norm regularization is closely related to Sparse Principal Component Analysis, which is a variant of PCA that aims to find a sparse principal components. It has real application in variable selection. \cite{jolliffe2003modified}, one of the first papers on sparse PCA, accomplished this goal through a formulation that maximized the quadratic form, incorporating a simplified constraint of the $\ell^1$ norm. In contrast, \cite{d2004direct} suggested estimating the sparse best-fit line by resolving a convex relaxation of a variance maximization issue. \cite{zou2006sparse}, on the other hand, proposed using a simplified version of elastic net regression on principal scores to achieve a sparse best-fit line.

\section{Estimating an $\ell^1$-Norm Regularized $\ell^1$-Norm Best-Fit Line}
In this section, we will extend the sorting method introduced in \cite{Brooks2017IdentifyingRanking} for estimating L1-norm best-fit lines to the setting where we add a penalty for sparsity. First, we introduce four sets of goal variables $\epsilon_{ij}^+,\epsilon_{ij}^-$ and $\zeta_j^+,\zeta_j^-$.  The optimization problem in $\eqref{formulation1}$ can be recast as the following constrained mathematical program.
    \begin{align}
        \label{formulation2}
         \displaystyle\min_{\genfrac{}{}{0pt}{2}{v\in \mathbb{R}^m,\alpha \in \mathbb{R}^n}{ \genfrac{}{}{0pt}{2}{\epsilon^+, \epsilon^- \in \mathbb{R}^{n\times m}_+,}{\zeta^+, \zeta^- \in \mathbb{R}^m_+}} } &   \sum_{i\in N}\sum_{j\in M} (\epsilon_{ij}^+ + \epsilon_{ij}^-) + \lambda\sum_{j \in M} (\zeta_j^+ + \zeta_j^-), 
    \end{align}
s.t.
    \begin{align*} 
         v_j\alpha_i+\epsilon_{ij}^+-\epsilon_{ij}^-&=x_{ij}, i\in N, j\in M, \\
         v_j + \zeta_j^+ - \zeta_j^-&=0, j \in M,\\
         \epsilon_{ij}^+, \epsilon_{ij}^-, \zeta_j^+, \zeta_j^-  & \geq 0, i \in N, j\in M. 
    \end{align*} 
\begin{prop}
The formulation \eqref{formulation2} is equivalent to \eqref{formulation1}.
\end{prop}  
\begin{proof}
The presence of absolute values in the objective function can be avoided by replacing each $x_{ij}-v_j\alpha_i$ with $\epsilon_{ij}^+-\epsilon_{ij}^-, \epsilon_{ij}^+,\epsilon_{ij}^-\geq 0$ and each $v_j$ with $\zeta_j^+-\zeta_j^-,\zeta_j^+,\zeta_j^-\geq 0$, and these become the constraints. The new objective function $\sum_{i=1}^n\sum_{j=1}^m|\epsilon_{ij}^+ - \epsilon_{ij}^-| + \lambda |\zeta_j^+ - \zeta_j^-|$ can be replaced with $\sum_{i=1}^n\sum_{j=1}^m(\epsilon_{ij}^+ + \epsilon_{ij}^-) + \lambda (\zeta_j^+ + \zeta_j^-)$. This linear program will have an optimal solution with at least one of the values in $\epsilon_{ij}^+, \epsilon_{ij}^-$ and $\zeta^+, \zeta^-$ is zero respectively. In that case, $x_{ij}-v_j\alpha_i=\epsilon_{ij}^+$, if $x_{ij}-v_j\alpha_i>0$, and $x_{ij}-v_j\alpha_i=-\epsilon_{ij}^-$, if $x_{ij}-v_j\alpha_i<0$. $v_j=\zeta_j^+$, if $v_j>0$, and $v_j=-\zeta_j^-$, if $v_j<0$. Any feasible solution for \eqref{formulation2} generates an objective function value which is the same as that of \eqref{formulation1} using the same values for $v$ and $\alpha$, and vice-versa. Therefore, an optimal solution to \eqref{formulation1} generates a feasible solution for \eqref{formulation2} and vice-versa.
\end{proof} 

An optimal solution to \eqref{formulation2} will be a vector $v^* \in \mathbb{R}^m$, along with scalars $\alpha_i^*$, $i\in N$.  For each point $i$ and feature $j$, the pair ($\epsilon_{ij}^{+*}$,$\epsilon_{ij}^{-*}$) reflects the distance along the unit direction $j$ between the point and its projection.  The pairs $(\zeta_j^{+*}, \zeta_j^{-*})$ provide the difference from zero for each coordinate of $v^*$.

The following proposition provides a foundation for the sorting method proposed by \cite{Brooks2017IdentifyingRanking}.

\begin{prop}
\label{p2}
{\bf \cite{Brooks2019ApproximatingL1}} Let $v\neq$ 0 be a given vector in $\mathbb{R}^m$. Then there is an $\ell^1$-norm projection from the point $x_i\in \mathbb{R}^m$ on the line defined by v that can be reached using at most $m-1$ unit directions.  Moreover, if the preserved coordinate is $\hat{\jmath}$ and $x_{i\hat{\jmath}} \neq 0$, then $v_{\hat{\jmath}} \neq 0$.
\end{prop}
\begin{proof}
A proof is in \cite{Brooks2019ApproximatingL1}.
\end{proof}
 
The idea is to impose the preservation of the same coordinate, $\hat{\jmath}$, in the projections of all points. \cite{Tsagkarakis2017OnRank,chier2017} also propose methods for subspace estimation based on the assumption that all points preserve the same unit directions.  In this work, we are focused on line fitting and we are adding a regularization term to promote sparsity.  

For a line, preservation of the same coordinate means that each point will use the same $m-1$ unit directions to project onto the line defined by $v$. 
 
By Proposition \ref{p2}, if $x_{\hat{\jmath}}\neq 0$, then $v_{\hat{\jmath}} \neq0$. Therefore, we can set $v_{\hat{\jmath}}=1$ and set $\alpha_i=x_{i\hat{\jmath}}$ to preserve $\hat{\jmath}$ without loss of generality for the error term, though the regularization term is affected.

The remaining components of $v$ can be found by solving an LP based on \eqref{formulation2} after replacing $\alpha_i$ with $x_{i\hat{\jmath}}$:
\begin{align}
        \label{formulation3}
        z_{\hat{\jmath}}(\lambda) =\displaystyle\min_{\genfrac{}{}{0pt}{2}{v\in \mathbb{R}^m,v_{\hat{\jmath}}=1 }{\genfrac{}{}{0pt}{2}{\epsilon^+, \epsilon^- \in \mathbb{R}^{n\times m},}{\zeta^+, \zeta^- \in \mathbb{R}^m}}}  & \sum_{i\in N}\sum_{j\in M} (\epsilon_{ij}^+ + \epsilon_{ij}^-) + \lambda \sum_{j\in M}(\zeta_j^+ + \zeta_j^-), 
    \end{align}
s.t.
    \begin{align*} 
        v_jx_{i\hat{\jmath}}+\epsilon_{ij}^+-\epsilon_{ij}^-&=x_{ij}, i\in N, j\in M;j\neq \hat{\jmath}, \\
        v_j + \zeta_j^+ - \zeta_j^-&=0, j\in M,\\
         \epsilon_{ij}^+, \epsilon_{ij}^-, \zeta_j^+, \zeta_j^-  & \geq 0, i \in N, j\in M.
    \end{align*}
Each of the $n$ data points generates $m-1$ constraints in this LP. 

Allowing each coordinate $j$ to serve as the preserved coordinate $\hat{\jmath}$ produces $m$ LPs.  By solving these $m$ LPs and selecting the vector $v$ from the solutions associated with the smallest values of the objective function $m$, we will have the $\ell^1$-norm regularized $\ell^1$-norm best-fit line under the assumption that all points project by preserving the same coordinate $\hat{\jmath}$ and $v_{\hat{\jmath}}=1$.  The following lemma describes how to generate solutions to the LPs by sorting several ratios.
\begin{lemma}
    \label{p4}
    For data $x_i \in \mathbb{R}^m$, $i\in N$, and for a $\lambda \in \mathbb{R}$, an optimal solution to \eqref{formulation3} can be constructed as follows.  If $x_{i\hat{\jmath}} = 0$ for all $i$, then set $v = 0$.  Otherwise, set $v_{\hat\jmath}=1$ and for each $j \neq \hat{\jmath}$,
    \begin{itemize}
        \item Take points $x_i$, $i\in N$ such that $x_{i\hat{\jmath}} \neq 0$ and sort the ratios $\displaystyle\frac{x_{ij}}{x_{i\hat{\jmath}}}$ in increasing order. 
        \item If there is a point $\tilde{\imath}$ where 
       \begin{equation} 
		    \label{solncond}
		    \left|\mbox{sgn}\left(\frac{x_{\tilde{\imath}j}}{x_{\tilde{\imath}\hat{\jmath}}}\right) \lambda +  \sum_{\genfrac{}{}{0pt}{2}{i \in N:}{i < \tilde{\imath}}} |x_{i\hat{\jmath}}| - \sum_{\genfrac{}{}{0pt}{2}{i \in N:}{i > \tilde{\imath}}} |x_{i\hat{\jmath}}| \right| \leq  |x_{\tilde{\imath}\hat{\jmath}}|,
		  \end{equation} 
    then set $v_j = \displaystyle\frac{x_{\tilde{\imath}j}}{x_{\tilde{\imath}\hat{\jmath}}}$.
               \item If no such $\tilde{\imath}$ exists, then set $v_j = 0$.
    \end{itemize}
\end{lemma}
\begin{proof} 
    The problem \eqref{formulation3} is separable into $m$ independent small sub-problems, one for each column $j$. 
    For a given $j$, the there is an LP of the form
    \begin{align}
        \label{l1reglp}
        \min_{\displaystyle\genfrac{}{}{0pt}{2}{v_j, \epsilon^+, \epsilon^-, \lambda}{ \zeta^+, \zeta^-}}  & \sum_{i\in N} (\epsilon_{ij}^+ + \epsilon_{ij}^-) + \lambda (\zeta_j^+ + \zeta_j^-), \\
        \mbox{s.t. } &         v_j x_{i\hat{\jmath}}+\epsilon_{ij}^+-\epsilon_{ij}^- =  x_{ij}, i\in N, \nonumber \\
        & v_j + \zeta_j^+ - \zeta_j^- = 0, \nonumber \\
        & \epsilon_{ij}^+, \epsilon_{ij}^-,\zeta_j^{+}, \zeta_j^- \geq 0, i \in N. \nonumber
    \end{align}

    We will show that the solution for $v_j$ stated in Lemma \ref{formulation1} is optimal by constructing a dual feasible solution that is complementary to the proposed primal feasible solution.  Suppose that the ratios $\frac{x_{ij}}{x_{i\hat{\jmath}}}$, $i\in N$, are sorted in increasing order.

    The dual linear program to \eqref{l1reglp} is
    \begin{align}
        \max_{\pi, \gamma} & \sum_{i \in N} \frac{x_{ij}}{x_{i\hat{\jmath}}} \pi_{i}, \\
        \mbox{s.t. } & \sum_{i \in N} \pi_{i} + \gamma = 0, \\
        & -|x_{i\hat{\jmath}}| \leq \pi_{i} \leq |x_{i\hat{\jmath}}|, i \in N, \label{pibound} \\
        & -\lambda \leq \gamma \leq \lambda.
    \end{align}
    Suppose there is an $\tilde{\imath}$ satisfying \eqref{solncond}.  Then let $\gamma = -\mbox{sgn} \left(\frac{x_{\tilde{\imath}j}}{x_{\tilde{\imath}{\hat{\jmath}}}}\right)\lambda$ and let
    \[
        \pi_i = \left\{\begin{array}{rl} |x_{i\hat{\jmath}}| & \mbox{if $i > \tilde{\imath}$},\\
                                - |x_{i\hat{\jmath}}| & \mbox{if $i < \tilde{\imath}$},\\
                                - \gamma - \sum_{i \neq \tilde{\imath}} \pi_i & \mbox{if $i = \tilde{\imath}$}.
        \end{array}\right.
    \]
    This solution satisfies complementary slackness.  To show that the solution is dual feasible, we need to show that $\pi_{\tilde{\imath}}$ satisfies the bounds in \eqref{pibound} (all other bounds and constraints are satisfied):  
        \begin{eqnarray}
            |\pi_{\tilde{\imath}}| & = & |-\gamma - \sum_{i \neq \tilde{\imath}} \pi_i|, \\
                                 & = & \left|\mbox{sgn} \left(\frac{x_{\tilde{\imath}j}}{x_{\tilde{\imath}{\hat{\jmath}}}}\right)\lambda  + \sum_{i: i < \tilde{\imath}} |x_{i\hat{\jmath}}| - \sum_{i: i > \tilde{\imath}} |x_{i\hat{\jmath}}| \right|,   \\ 
                                 & \leq & |x_{\tilde{\imath}\hat{\jmath}}|.
        \end{eqnarray}
  The inequality is due to \eqref{solncond}.  There is a complementary dual feasible solution, so the proposed solution must be optimal.

Now suppose that there is no $\tilde{\imath}$ satisfying \eqref{solncond}.  
Note that if \eqref{solncond} is satisfied for some $\tilde{\imath}$ with $\mbox{sgn}\left(\frac{x_{\tilde{\imath}j}}{x_{\tilde{\imath}\hat{\jmath}}}\right) = +$, then
\begin{align}
  \lambda & \geq \sum_{i: i > \tilde{\imath}} |x_{i\hat{\jmath}}| - \sum_{i: i \leq \tilde{\imath}} |x_{i \hat{\jmath}}|, \label{lambdalow}\\
  \lambda & \leq \sum_{i: i \geq \tilde{\imath}} |x_{i\hat{\jmath}}| - \sum_{i: i < \tilde{\imath}} |x_{i \hat{\jmath}}|. \label{lambdahigh}
\end{align}
If \eqref{solncond} is violated for each $\tilde{\imath}$, then for each possible $\tilde{\imath}$ either the lower bound \eqref{lambdalow} or the upper bound \eqref{lambdahigh} for $\lambda$ is violated.  If for a given $\tilde{\imath}$, the lower bound \eqref{lambdalow} is violated, then $\lambda  <  \sum_{i: i > \tilde{\imath}} |x_{i\hat{\jmath}}| - \sum_{i: i \leq \tilde{\imath}} |x_{i \hat{\jmath}}|$.  This implies that the upper bound \eqref{lambdahigh} is satisfied.  If we now consider point $\tilde{\imath} + 1$, then the upper bound is the same as the lower bound for $\tilde{\imath}$ and is therefore satisfied.  So $\lambda$ must violate the lower bound for $\tilde{\imath} + 1$, and we can consider $\tilde{\imath} + 2$ and so on.  Then lower bound is violated for all points with $\mbox{sgn}\left(\frac{x_{\tilde{\imath}j}}{x_{\tilde{\imath}\hat{\jmath}}}\right) = +$, in particular the largest, and so $\lambda < 0$, contradicting the choice of $\lambda$.  A symmetric argument  holds for $\tilde{\imath}$ with $\mbox{sgn}\left(\frac{x_{\tilde{\imath}j}}{x_{\tilde{\imath}\hat{\jmath}}}\right) = -$.  
Therefore, $\lambda > \sum_{i: i \geq \tilde{\imath}}|x_{i\hat{\jmath}}| - \sum_{i: i < \tilde{\imath}} |x_{i\hat{\jmath}}|$, for every $\tilde{\imath}$ with $\mbox{sgn}\left(\frac{x_{\tilde{\imath}j}}{x_{\tilde{\imath}\hat{\jmath}}}\right) = +$ and $\lambda > \sum_{i: i \leq \tilde{\imath}}|x_{i\hat{\jmath}}| - \sum_{i: i > \tilde{\imath}} |x_{i\hat{\jmath}}|$ for every $\tilde{\imath}$ with $\mbox{sgn}\left(\frac{x_{\tilde{\imath}j}}{x_{\tilde{\imath}\hat{\jmath}}}\right) = -$.  In particular, 
\begin{equation}
  \lambda > 
  \left|\sum_{i: \frac{x_{ij}}{x_{i\hat{\jmath}}} < 0} |x_{i\hat{\jmath}}|
  - \sum_{i: \frac{x_{ij}}{x_{i\hat{\jmath}}} > 0} |x_{i\hat{\jmath}}| \right|.
\end{equation}
A dual feasible and complementary solution is to set 
    \[
      \pi_i = \left\{\begin{array}{rl} |x_{i\hat{\jmath}}| & \mbox{if} \frac{x_{ij}}{x_{i\hat{\jmath}}} > 0,\\
	- |x_{i\hat{\jmath}}| & \mbox{if} \frac{x_{ij}}{x_{i\hat{\jmath}}} < 0,
        \end{array}\right.
    \]
and $\gamma = \left|\sum_{i: \frac{x_{ij}}{x_{i\hat{\jmath}}} < 0} |x_{i\hat{\jmath}}|
  - \sum_{i: \frac{x_{ij}}{x_{i\hat{\jmath}}} > 0} |x_{i\hat{\jmath}}| \right|$.  Note that $|\gamma| < \lambda$ by the development above, so the solution is dual feasible and therefore optimal.  
 
 \end{proof}
 
Given a penalty $\lambda$, an optimal solution to \eqref{formulation3} with the preservation of one coordinate $\hat{\jmath}$ and $v_{\hat{\jmath}} = 1$ requires the sorting of $(m-1)$ lists of ratios according to Lemma \ref{p4}. The process is repeated for each choice of $\hat{\jmath}$ and the solution with the smallest objective function value is retained. Therefore, for a penalty $\lambda$, the proposed method requires sorting $m(m-1)$ lists of ratios in total, each costing $(n\log n)$ running time.  This motivates an $O{(m^2n \log n)}$ algorithm for estimating $v$ formalized in Algorithm 1 below. 
\begin{prop}
For a given $\lambda$ and data $x_i \in \mathbb{R}^m$, $i\in N$, Algorithm 1 finds an optimal solution to \eqref{formulation3}. 
\end{prop}   
\begin{proof}
For each fixed coordinate, Algorithm 1 finds an optimal solution according to Lemma 1.  From among those solutions, Algorithm 1 picks the one with the smallest combination of error plus regularization term.
\end{proof}

\begin{algorithm}[!htbp]
 \caption{Estimating an $\ell^1$-norm regularized $\ell^1$-norm best-fit line $v^*$ for given $\lambda$.}
 \algsetup{indent=2em} 
 \begin{algorithmic}[1]
 \renewcommand{\algorithmicrequire}{\textbf{Input:}}
 \renewcommand{\algorithmicensure}{\textbf{Output:}} 
 \newcommand{\algorithmicbreak}{\textbf{break}}
 \newcommand{\BREAK}{\STATE \algorithmicbreak}
 \REQUIRE $x_i \in \mathbb{R}^{m}$ for $i = 1,\dotso,n$. $\lambda$.
 \ENSURE $v^*$
  \STATE{Set $z^*=\infty$}
   \FOR {$\hat{\jmath}\in M$}
    \STATE{Set $v_{\hat{\jmath}} = 1$.} 
  \FOR{$j\in M: j \neq\hat{\jmath}$}
  \STATE {Set $v_j = 0$.}
  \STATE{Sort $\left\{\frac{x_{ij}}{x_{i\hat{\jmath}}}: i\in N, x_{i\hat{\jmath}} \neq 0\right\}$.}  
    \FOR{$\tilde{\imath} \in N: x_{\tilde{\imath}\hat{\jmath}} \neq 0$} 
      \IF{$\mbox{sgn}\left(\frac{x_{ij}}{x_{i\hat{\jmath}}}\right)\lambda\in(\sum\limits_{i: i > \tilde{\imath}} |x_{i\hat{\jmath}}|-\sum\limits_{i: i \leq \tilde{\imath}} |x_{i\hat{\jmath}}|,\sum\limits_{i: i \geq \tilde{\imath}} |x_{i\hat{\jmath}}|-\sum\limits_{i: i < \tilde{\imath}} |x_{i\hat{\jmath}}|]$}  
          \STATE Set $v_j=\frac{x_{\tilde{\imath}j}}{x_{\tilde{\imath}\hat{\jmath}}}$. 
       \ENDIF
      \ENDFOR 
  \ENDFOR
  \STATE {set $z=\sum\limits_{i\in N}\sum\limits_{j\in M}|x_{ij}-v_j x_{i\hat{\jmath}}|+\lambda\sum\limits_{j\in M}|v_{j}|$}
  \IF{$z<z^*$}
  \STATE{Set $z^* = z$, $v^{*}=v$}
  \ENDIF
  \ENDFOR
  \RETURN {$v^*$}
 \end{algorithmic} 
 \label{alg1}
 \end{algorithm} 
 
Algorithm 1 finds the best solution that preserves each coordinate $\hat{\jmath}$ and $v_{\hat{\jmath}}$=1 for a given value of $\lambda$. Algorithm 2 seeks the intervals constructed by successive breakpoints - values for $\lambda$ at which the solution is going to change and the conditions of Lemma 1 are satisfied. Algorithm 2 does not determine which coordinate $\hat{\jmath}$ is best to preserve for each interval.
 Algorithm 3 iterates through each interval for $\lambda$ from Algorithm 2 and finds the intervals where preserving $\hat{\jmath}$ minimizes the objective function value.

\begin{algorithm} [!htbp] 
 \caption{Find all major breakpoints.}
 \begin{algorithmic}[1]
 \renewcommand{\algorithmicrequire}{\textbf{Input:}}
 \renewcommand{\algorithmicensure}{\textbf{Output:}}
 \REQUIRE $x_i \in \mathbb{R}^{m}$ for $i = 1,\dotso,n$.
 \ENSURE Ordered breakpoints for the penalty $\Lambda$ and solutions $v^{\hat{\jmath}}(\lambda)$ for each choice of preserved coordinate $\hat{\jmath}$, and each $\lambda \in \Lambda$.
 \STATE Set $\Lambda = \{0, \infty\}$.
  \FOR {$\hat{\jmath}\in M$}
    \STATE{Set $v^{\hat{\jmath}}_{\hat{\jmath}} = 1$.}
  \FOR{$j\in M: j \neq\hat{\jmath}$}
  \STATE{Set $\lambda^{\max} = 0$.}
  \STATE{Sort $\left\{\frac{x_{ij}}{x_{i\hat{\jmath}}}: i\in N, x_{i\hat{\jmath}} \neq 0\right\}$.}
    \FOR{$\tilde{\imath} \in N: x_{\tilde{\imath}\hat{\jmath}} \neq 0$}
       \STATE Set {$\lambda = \mbox{sgn}\left(\frac{x_{\tilde{\imath}j}}{x_{\tilde{\imath}\hat{\jmath}}}\right) \left(\sum\limits_{i: i > \tilde{\imath}} |x_{i\hat{\jmath}}|-\sum\limits_{i: i < \tilde{\imath}} |x_{i\hat{\jmath}}|\right) - |x_{\tilde{\imath}\hat{\jmath}}|$}
       \IF{$\lambda + 2|x_{\tilde{\imath}\hat{\jmath}}| > 0$,}  
         \STATE Set $\Lambda = \Lambda \cup \max\{0,\lambda\}$.
         \STATE Set $v^{\hat{\jmath}}_j(\max\{0,\lambda\}) = 
         \frac{x_{\tilde{\imath}j}}{x_{\tilde{\imath}\hat{\jmath}}}$.
       \ENDIF
       \IF{$\lambda +2|x_{\tilde{\imath}\hat{\jmath}}| > \lambda^{\max}$,}
         \STATE{Set $\lambda^{\max} = \lambda + 2|x_{\tilde{\imath}\hat{\jmath}}|$.}
       \ENDIF
     \ENDFOR
     \STATE{Set $\Lambda = \Lambda \cup \{\lambda^{\max}\}$.}
     \STATE{Set $v_j^{\hat{\jmath}}(\lambda^{\max}) = 0$.}
  \ENDFOR
  \ENDFOR 
  \STATE Sort $\Lambda$.
  \RETURN $\Lambda$, $\{v_j^{\hat{\jmath}}(\lambda): j \in M, \hat{\jmath} \in M, \lambda \in \Lambda\}$
 \end{algorithmic} 
 \label{alg2}
 \end{algorithm}  
 
 It is necessary to ''merge'' the intervals for each possible preserved coordinate $\hat{\jmath}$ and determine when the preservation of each coordinate results in the lowest value of the objective function.  Therefore, we need Algorithm 3 to check each consecutive interval for $\lambda$ from Algorithm 2 to determine if changing the preserved coordinate $\hat{\jmath}$ can reduce the objective function value, which may result in new breakpoints that were not discovered using Algorithm 2. 
\begin{algorithm} [!htbp] 
 \caption{Solution Path for $\ell^1$-norm Regularized $\ell^1$-norm best-fit line}
 \begin{algorithmic}[1]
 \renewcommand{\algorithmicrequire}{\textbf{Input:}}
 \renewcommand{\algorithmicensure}{\textbf{Output:}} 
 \REQUIRE A ordered set of breakpoints for the penalty  $(\lambda^k: k=1,\ldots, K)$ and solutions $v^{\hat{\jmath}}(\lambda^k)$ for each choice of preserved coordinate $\hat{\jmath}$. 
 \ENSURE Breakpoints for the penalty $\Lambda$ and solutions $v^*(\lambda)$ for each $\lambda \in \Lambda$. 
 \STATE{Set $\Lambda = \emptyset.$}
  \FOR{$k = 1, \ldots, K-1$}
    \FOR{$\hat{\jmath} \in M$}
      \STATE{Set $z^{\hat{\jmath}}(\lambda^k) = \sum_{i \in N} \|x_i - v^{\hat{\jmath}}x_{i\hat{\jmath}}\|_1 + \lambda^k \|v^{\hat{\jmath}}(\lambda^k)\|_1$}
    \ENDFOR
    \FOR{$j \in M$} 
      \STATE $\beta_L = 
      \begin{aligned}[t]
         &\left\{ \max \frac{z^{j}(\lambda^k) - z^{\hat{\jmath}}(\lambda^k)}{\|v^{\hat{\jmath}}(\lambda^k)\|_1 - \|v^j(\lambda^k)\|_1}:\right. \\
         &\left. \hat{\jmath} \vphantom{\bigg\{} \in M, \|v^j(\lambda^k)\|_1 < \|v^{\hat{\jmath}}(\lambda^k)\|_1  \right\}
      \end{aligned}$
      \STATE $\beta_U = 
      \begin{aligned}[t]
        &\left\{ \{\min \frac{z^{\hat{\jmath}}(\lambda^k) - 
        z^{j}(\lambda^k)}{\|v^{j}(\lambda^k)\|_1 - 
        \|v^{\hat{\jmath}}(\lambda^k)\|_1}:\right.  \\
        &\left. \hat{\jmath} \vphantom{\bigg\{} \in M, \|v^j(\lambda^k)\|_1 > \|v^{\hat{\jmath}}(\lambda^k)\|_1 \right\} 
      \end{aligned}$
      \IF{$\left|\left\{\hat{\jmath}: z^j(\lambda^k) > z^{\hat{\jmath}}(\lambda^k), \|v^j(\lambda^k)\|_1 = \|v^{\hat{\jmath}}(\lambda^k)\|_1\right\}\right| = 0,$}
        \IF{$0< \beta_L < \beta_U$ and $\lambda^k + \beta_L \leq \lambda^{k+1}$, }
        \STATE{Set $\Lambda = \Lambda \cup \{\lambda^k + \beta_L\}$}
        \STATE{Set $v^*(\lambda^k + \beta_L) = v^j(\lambda^k)$}
        \ELSIF{$\beta_L \leq 0 < \beta_U$,}
          \STATE{Set $\Lambda = \Lambda \cup \{\lambda^k\}$}
          \STATE{Set $v^*(\lambda^k) = v^j(\lambda^k)$}
        \ENDIF
      \ENDIF
    \ENDFOR
  \ENDFOR
  \RETURN $\Lambda$, $\{v^*(\lambda): \lambda \in \Lambda\}$
 \end{algorithmic} 
 \label{alg3}
 \end{algorithm}
 
\begin{prop}
\label{p5} 
For data $x_i \in \mathbb{R}^m$, $i\in N$, Algorithms 2 and 3 generate the entire solution path for \eqref{formulation3} across all possible values of $\lambda$ under the assumption that all points are projected, preserving the same unit direction and $v_{\hat{\jmath}}$=1 for the preserved direction $\hat{\jmath}$.  
\end{prop}   

 \section{Experiments with Synthetic Data}
In this section, we shall first shift our attention by analyzing a toy sample, trying to understand the complete solution path in terms of breakpoints and coordinate preservation, that is, how breakpoints affect solutions by changing preserved coordinates $\hat{\jmath}$. Next, we conducted simulation studies to evaluate the performance of Algorithm \ref{alg1} against some classic low-rank approximation algorithms.  Error is measured by the discordance between the vector $v$ defining the ``true'' line and the vector $v^*$ derived by Algorithm \ref{alg1} or a competing method.  Sparsity is measured by the $\ell^0$ norm of the solution vector $v$ which is the number of non-zero coordinates.

\subsection{A Toy Example}
\label{sec:toy}
Let us first consider five points $(4,-2,3,-6)^T$, $(-3,4,2,-1)^T$, $(2,3,-3,-2)^T$, $(-3,4,2,3)^T$, $(5,3,2,-1)^T$. Algorithm 2 generates the following breakpoints for $\lambda$: $\{1, 3\}$ for $\hat{\jmath}=1$, $\{4, 6\}$ for $\hat{\jmath}=2$, $\{0, 2\}$ for $\hat{\jmath}=3$ and $\{3, 5, 11\}$ for $\hat{\jmath}=4$. The collection of breakpoints for $\lambda$ is $\{0,1,2,3,4,5,6,11\}$. We now illustrate that the optimal solution (under the assumption that all points preserve the same coordinate) might change due to the existence of additional breakpoints between successive breakpoints generated from Algorithm 2. Algorithm 3 iterates by preserving $j=1,2,3$ to find the lowest objective value over the interval $(3.5,4]$ for $\lambda$, giving rise to an additional breakpoint $3.5$. The value of the objective function comprises the error term ($\sum_{i=1}^n\|x_{i}-vx_{i\hat{\jmath}}\|_1$) and the penalty term ($\|v\|_1$), both fixed over each interval for each coordinate $j$. Algorithm 2 finds all possible breakpoints without filtering comparable larger objective function values, which is assessed in Algorithm 3. In other words, Algorithm 3 further narrows the breakpoint intervals of Algorithm 2 by evaluating $m$ objective function values. The complete solution path is summarized in Table \ref{tab:Of}. 
\begin{table}[!t]  
\centering
  \caption{Solution Path for Toy Example. The Best-Fit Line Is Fixed within Each of the Four Intervals for $\lambda$.}  
  \begin{tabular}{r|r|r} 
     $\lambda$&$z^*(\lambda)$& $v^*(\lambda)$\\
    \hline 
   (0.0, 3.0)   & (34.5, 42.0) & (-0.7,0.3,-0.5,1.0)\\ 
   (3.0, 3.5)   & (42.0, 42.9) & ($-\frac{2}{3}$,$\frac{1}{3}$,0.0,1.0)\\  
   (3.5, 11) & (42.9, 52.0) &  (1.0,0.0,0.0,-0.2)  \\ 
        (11, $\infty$) & (52.0, $\infty$) & (1.0,0.0,0.0,0.0) \\
\end{tabular} 
\label{tab:Of}
\end{table}

\begin{figure} [!t] 
\centering 
\begin{tikzpicture}  
\begin{axis}[%
    width=0.50\columnwidth,
    scale only axis,
    compat=newest,
    ylabel style = {align=left},
    axis line style=thick,
    inner axis line style={=>},
    xlabel={$\lambda$},  
    ylabel={$z_{\hat{\jmath}}$}, 
      x label style={anchor=west},
        y label style={anchor=south},
    ymin=34,ymax=58,
    legend pos=north west,
    xmin=0,xmax=12,
    xticklabels={},
    extra x ticks = {1,2,3,4,5,6,11},
    ticklabel style = {font=\scriptsize},
    legend style={font=\scriptsize},
    legend entries={$\hat{\jmath}=1$,
    $\hat{\jmath}=2$,
    $\hat{\jmath}=3$,
    $\hat{\jmath}=4$
    }
]
\addplot[line width=0mm ]{0};
\addplot[line width=0mm,cyan]{0}; 
\addplot[line width=0mm,blue]{0};  
\addplot[line width=0mm,color=red]{0};

\addplot[line width=0.3mm,domain=0:1 ]{36.1+2.9*x};
\addplot[line width=0.3mm,domain=1:3 ]{37.3+1.7*x}; 
\addplot[line width=0.3mm,domain=3:11 ]{38.8+1.2*x};  
\addplot[line width=0.3mm,domain=11:12 ]{41+1*x};

\addplot[line width=0.3mm,domain=0:4, cyan]{35+2.5*x};
\addplot[line width=0.3mm,domain=4:6, cyan,]{39+1.5*x};  
\addplot[line width=0.3mm,domain=6:12, cyan,]{42+1*x};  

\addplot[line width=0.3mm,domain=0:2, blue]{43.7+2.17*x};
\addplot[line width=0.3mm,domain=2:12, blue]{46+x}; 

\addplot[line width=0.3mm,domain=0:3,red,]{34.5+2.5*x};
\addplot[line width=0.3mm,domain=3:5,red,]{36+2*x};    
\addplot[line width=0.3mm,domain=5:11,red,]{37.66666+1.66666*x};
\addplot[line width=0.3mm,domain=11:12,red]{45+x};     

\draw[dashed,gray] (1,0) -- (1,39);
\draw[dashed,gray] (2,0) -- (2,48);
\draw[dashed,gray] (3,0) -- (3,42.4);
\draw[dashed,gray] (4,0) -- (4,45);
\draw[dashed,gray] (5,0) -- (5,46.2);
\draw[dashed,gray] (6,0) -- (6,48);
\draw[dashed,gray] (11,0) -- (11,52);
\draw[gray] (3.5,0) -- (3.5,42.9);
\end{axis}
\end{tikzpicture}
\caption{Schematic illustration of breakpoints. Each color depicts the objective function value when preserving a coordinate $\hat{\jmath}$, $z_{\hat{\jmath}}$, as a function of the penalty parameter $\lambda$.} 
\label{fig:ex-b}
\end{figure} 

The objective function $z_{\hat{\jmath}}$ is a linear function with respect to $\lambda$ over a certain interval for each preserved $\hat{\jmath}$. The intercept is $\sum_{i=1}^n \|x_{i}-vx_{i\hat{\jmath}}\|_1$, and $\|v\|_1$ is the slope. Figure \ref{fig:ex-b} shows that the smallest $z_{\hat{\jmath}}$ can be achieved by preserving direction $\hat{\jmath}=4$ for $\lambda$ in the interval $[0,3.5]$, depicted by the red line segment, and direction $\hat{\jmath}=1$ for $\lambda$ in the interval $[0,\infty)$, depicted by the black line segment.  Algorithm 2 finds all 7 breakpoints $\lambda$s at which the slope of the lines in same color changes. Algorithm 3 finds an additional breakpoint in the interval (3,4], depicted by the gray vertical line, at which value two line segments intersect. At this value, the preserved direction changes from $\hat{\jmath}=4$ to $\hat{\jmath}=1$ as $\lambda$ increases.  

This result is consistent with that of Table \ref{tab:Of}, where the first two solutions preserve $\hat{\jmath}=4$ and the last two solutions preserve $\hat{\jmath}=1$. The solutions for each choice of $\lambda$ are summarized in Table \ref{tab:Of}.

\begin{align}
    x_{\hat{\jmath}}v^*&=v'\alpha' \nonumber \\
    &= \frac{v^*}{\|v^*\|}\alpha'\nonumber \\
    \alpha' &= x_{\hat{\jmath}}\|v^*\|\nonumber
\end{align}

\begin{remark}
For any solution to \eqref{formulation1}, a solution with a better objective function can be obtained by dividing $v$ by a large constant.  Therefore, $v$ needs to be normalized in magnitude in some way.  We will do so by requiring $v_{\hat{\jmath}} = 1$.   
\end{remark} Successive components can be computed by applying Algorithm \ref{alg1} to the data projected in the null space of the subspace defined by the components derived thus far. 
 
\subsection{Comparison to Competing Methods}

In this section, we demonstrate the effectiveness of Algorithm \ref{alg1} compared to principal component pursuit \cite{Candes2011RobustAnalysis} (PCP), augmented Lagrange multiplier\cite{tang2011robust} (LRBS) and alternating direction multiplier method \cite{tao2011recovering} (ADM). For each configuration, a total of 10 replications are performed. For each replication, we create datasets with $n$ observations $\mathbb{R}^m$ including $nC$ outliers. Therefore, the values of $n$ and $m$ are the number of rows and the number of columns of input data. The values of $nC$ and $mC$ are the number of rows and columns contaminated. Each element of the ``true'' $v$ is sampled from a Uniform distribution (-1,1) and $v$ is normalized.  To locate the projections of points on $v$, $\alpha_i$ is sampled from a Uniform distribution (-100,100) for $i=1,\ldots,n$. Synthetic data are generated by $\alpha_i v^T + \epsilon_i$, for $i=1,\ldots, n$, where $\epsilon_i$ is the noise sampled from a Laplace (0, 1) distribution with the probability density function $f(\epsilon_i|0,1)=0.5e^{-|\epsilon_i|}$. There is a link between the median and the Laplace distribution in the sense that the maximum likelihood estimator of location parameter for a list of independent and identically distributed samples following the Laplace distribution is the sample median \cite{li2004maximum}. Clustered outlier observations are created by sampling by first establishing an outlier center.  The first $mC$ coordinates of the center are sampled from a Uniform (100,150) distribution and the remainder are 0.   Outliers are created by adding noise 
 to the outlier center sampled from a Laplace (0,0.1) distribution.  All free parameters for PCP, LRBS, and ADM are set by default. The default value $\lambda$ for Algorithm \ref{alg1} is chosen as the average value of all breakpoints.   
 
\begin{table}[!t] 
\centering
\caption{Average Discordance Across 10 Replications.  The Subscripts Are the Standard Deviations. Values Less than 0.001 are denoted by $-$.}
\label{tab:discordance}
\begin{tabular}{cccc|ccccc}  
     $n$&$m$&$nC$&$mC$&PCA&PCP&LRBS&ASALM&Algorithm \ref{alg1}\\ 
     \hline
    1000&100&0&0&$-_{-}$&$-_{-}$&$-_{-}$&$-_{-}$&$-_{-}$\\  
  1000&100&100&5&$0.9_{.05}$&$-_{-}$&$-_{-}$&$0.9_{.13}$&$-_{-}$ \\  
  10000&100&0&0&$-_{-}$&$-_{-}$&$-_{-}$&$-_{-}$&$-_{-}$\\   
10000&100&1000&5&$0.8_{0.1}$&$-_{-}$&$0.87_{.09}$&$0.84_{.11}$&$-_{-}$ \\
   1000&1000&0&0&$-_{-}$&$-_{-}$&$-_{-}$&$-_{-}$&$-_{-}$\\  
 1000&1000&100&5&$0.9_{.04}$&$-_{-}$&$0.96_{.05}$&$0.9_{.08}$&$-_{-}$ \\ 
   1000&2000&0&0&$-_{-}$&$-_{-}$&$-_{-}$&$-_{-}$&$-_{-}$ \\ 
1000&2000&100&5&$0.9_{.06}$&$-_{-}$&$0.99_{.04}$&$0.98_{.08}$&$-_{-}$ \\ 
   5000&2000&0&0&$-_{-}$&$-_{-}$&$-_{-}$&$-_{-}$&$-_{-}$  \\ 
5000&2000&1000&5&$0.9_{.02}$&$-_{-}$&$0.98_{.02}$&$0.97_{.02}$&$-_{-}$ 
\end{tabular}

\end{table}

As can be seen in Table \ref{tab:discordance}, Algorithm \ref{alg1} and PCP produce accurate estimation in terms of low discordance over all configurations. (The cosine of the angle between two unit vectors $v_{est}$ and $v_{true}$ is equal to their dot product. Therefore, a smaller discordance implies a smaller angle between two unit vectors.) However, the precision of PCA, LRBS and ADM decreases significantly compared to that of Algorithm \ref{alg1}, when contamination is introduced into the data. In terms of $\ell^0$  in Table \ref{tab:nonzero}, the solutions of PCA and PCP do not exhibit any sparsity, whereas Algorithm \ref{alg1} produces more sparser solution for a given $\lambda$ without sacrificing much precision in presence of outliers. We will explore the effect of $\lambda$ on sparsity in later experiment.  ADM and LRBS produce solutions with some sparsity but with a large discordance. 
\begin{table} [!t]
\centering
\caption{The Average $\ell^0$ Sparsity Divided by $m$ Across 10 Replications. 
 The Subscripts Are the Standard Deviations.}
\begin{tabular}{cccc|ccccc}  
     $n$&$m$&$nC$&$mC$&PCA&PCP&LRBS&ASALM&Algorithm \ref{alg1}\\ 
    \hline 
    1000&100&0&0&$100_{0.0}$&$100_{0.0}$&$100_{0.0}$&$100_{0.0}$&$96.8_{0.8}$\\   
  1000&100&100&5&$100_{0.0}$&$100_{0.0}$&$100_{0.0}$&$99.8_{0.4}$&$97.2_{0.8}$ \\ 
     10000&100&0&0&$100_{0.0}$&$100_{0.0}$&$100_{0.0}$&$100_{0.0}$&$97.2_{2.4}$\\   
10000&100&1000&5&$100_{0.0}$&$100_{0.0}$&$99.9_{0.3}$&$99.9_{0.3}$&$96.2_{0.8}$ \\ 
   1000&1000&0&0&$100_{0.0}$&$100_{0.0}$&$100_{0.0}$&$100_{0.0}$&$98.3_{0.8}$\\   
 1000&1000&100&5&$100_{0.0}$&$100_{0.0}$&$98.5_{1.3}$&$98.8_{1.6}$&$91.7_{0.9}$ \\ 
   1000&2000&0&0&$100_{0.0}$&$100_{0.0}$&$100_{0.0}$&$100_{0.0}$&$90.2_{0.8}$ \\ 
1000&2000&100&5&$100_{0.0}$&$100_{0.0}$&$99.1_{1.4}$&$98.4_{1.9}$&$91.4_{1.0}$ \\ 
   5000&2000&0&0&$100_{0.0}$&$100_{0.0}$&$100_{0.0}$&$100_{0.0}$&$90.2_{0.7}$  \\ 
5000&2000&1000&5&$100_{0.0}$&$100_{0.0}$&$94.8_{3.7}$&$90.6_{12.0}$&$89.9_{0.6}$  
\end{tabular} 
\label{tab:nonzero}
\end{table}

\subsection{Effect of Varying $m$ and $n$ on Error and Sparsity} 
Positive regularization terms $\lambda$ are used to control the level of sparsity in the solution. And there shall be a maximum $\lambda$ beyond which all elements are 0 except for $v_{\hat{\jmath}}$ which is fixed at 1. In Figure \ref{fig:ll1}(a), the effects of varying $\lambda$ are illustrated in two settings, (a) varying the number of columns $m$ for a fixed number of rows $n$ and (b) varying the number of rows $n$ for a fixed number of columns $m$.  As can be seen, with the same number of rows, Algorithm \ref{alg1} is more sensitive to data with a larger number of columns in terms of the increasing rate of discordance and the decreasing rate of $\ell^0$. With the same number of columns, Algorithm \ref{alg1} is more sensitive to data with a smaller number of rows. It also shows that a $\lambda$ less than the intersection point (between the discordance curve, and the $\ell^0$-norm curve) will lead to a solution with discordance below 0.4 and sparsity below 40\%. Figure \ref{fig:ll1}(b) shows that solutions with similar properties can be obtained when varying $n$. Furthermore, Figures \ref{fig:ll1}(b) illustrates that Algorithm \ref{alg1} is more sensitive to $\lambda$ when working with data with smaller $m$. 

\begin{figure} [!t]  
\centering 
\begin{tikzpicture}[scale=1]
\pgfplotsset{ 
    xmin=0, xmax=8500
}
\begin{axis}[
    xlabel={$\lambda$},
    ylabel={Discordance/$\ell^0$},  
    legend style={draw=none},
    ymin=0, ymax=1, 
    legend style={at={(0.45,0.6)},anchor=north east},
    ymajorgrids=true,
    grid style=dashed, 
] 
\addplot[color=blue,mark=o,]
    table [x=lambda, y=1000d,col sep=comma] {csv/lbehavior_fixedrowsat3000.csv};
    \legend{3000x1000}
\addplot[color=red,mark=square,]
    table [x=lambda, y=2000d,col sep=comma] {csv/lbehavior_fixedrowsat3000.csv};
    \addlegendentry{3000x2000}
\addplot[color=green,mark=triangle,]
    table [x=lambda, y=3000d,col sep=comma] {csv/lbehavior_fixedrowsat3000.csv};
    \addlegendentry{3000x3000}
\addplot[color=blue,mark=o,dashed,mark options={solid}]
    table [x=lambda, y=1000l0,col sep=comma] {csv/lbehavior_fixedrowsat3000.csv};  
\addplot[color=red, dashed,mark=square,mark options={solid}]
    table [x=lambda, y=2000l0,col sep=comma] {csv/lbehavior_fixedrowsat3000.csv};  
\addplot[color=green,dashed,mark=triangle,mark options={solid}]
    table [x=lambda, y=3000l0,col sep=comma] {csv/lbehavior_fixedrowsat3000.csv};  
\end{axis}  
\end{tikzpicture}
\begin{tikzpicture}[scale=1]
\pgfplotsset{ 
    xmin=0, xmax=5500
}
\begin{axis}[
    xlabel={$\lambda$},  
    legend style={draw=none},
    ymin=0, ymax=1, 
    legend style={at={(0.45,0.6)},anchor=north east},
    ymajorgrids=true,
    grid style=dashed,
     ylabel={Discordance/$\ell^0$}, 
    ymin=0.0003, ymax=1,  
    yticklabel pos=left, 
] 
\addplot[color=blue,mark=o,]
    table [x=lambda, y=1000d,col sep=comma] {csv/lbehavior_fixedcolsat1000.csv};
    \legend{1000x1000}
\addplot[color=red,mark=square,]
    table [x=lambda, y=2000d,col sep=comma] {csv/lbehavior_fixedcolsat1000.csv};
    \addlegendentry{2000x1000}
\addplot[color=green,mark=triangle,]
    table [x=lambda, y=3000d,col sep=comma] {csv/lbehavior_fixedcolsat1000.csv};
    \addlegendentry{3000x1000}
\addplot[color=blue,mark=o,dashed,mark options={solid}]
    table [x=lambda, y=1000l0,col sep=comma] {csv/lbehavior_fixedcolsat1000.csv};  
\addplot[color=red, dashed,mark=square,mark options={solid}]
    table [x=lambda, y=2000l0,col sep=comma] {csv/lbehavior_fixedcolsat1000.csv};  
\addplot[color=green,dashed,mark=triangle,mark options={solid}]
    table [x=lambda, y=3000l0,col sep=comma] {csv/lbehavior_fixedcolsat1000.csv};  
\end{axis} 
\end{tikzpicture} 
\caption{Effect on discordance (solid lines) and sparsity (dashed lines) when $\lambda$ is varied for left datasets with different $m$ and right datasets with different $n$.  Sparsity is measured as a percentage of $m$.}
\label{fig:ll1}
\end{figure}

\section{Application to Human Microbiome Project Data}

This section demonstrates the efficacy of our solution sparsity in feature selection. In the context of clustering, feature selection searches for a small subset of features in the original collection as discriminators for the clustering task. Feature selection can increase the readability and interpretability of the model, which is of importance in many practical applications, such as finding genes most relevant to a specific disease and prognostic factors significant to treatment outcome prediction.

This experiment was carried out with data from the Human Microbiome Project \cite{peterson2009nih}. Data consist of relative abundances for 2,121 human tissue samples collected from 3 different source body sites (skin, gut, and oral). A total of 320 genera for each tissue sample is first fed into Algorithm 1 using different $\lambda$ values to get a subset of genera able to discriminate 3 body sites. A hierarchical clustering method with complete merging using Manhattan distance was carried on each subset. The performance of the clustering is measured by the purity defined as $\frac{1}{2355}\sum_{i=1}^{r}\max_{j=1}^k{c_{ij}}$ of $k=3$ body sites and $r=3$ clusters, where $c_{ij}$ represents the size of the group with which each cluster shares the most samples. 

\begin{table}[]
\centering
\caption{Core genera identified by Algorithm 1 with $\lambda$=224. Genus in bold is the single remaining discriminator. }
\label{tab:hmp}
\begin{tabular}{l|l|l|l}
\toprule
Actinomyces       & Parabacteroides & Streptococcus    & Subdoligranulum  \\
Corynebacterium   & Porphyromonas   & Eubacterium      & Veillonella      \\    
Rothia            & Alloprevotella  & Blautia          & Fusobacterium    \\  
Propionibacterium & Capnocytophaga  & Roseburia        & Leptotrichia     \\    
Atopobium         & Gemella         & Oscillibacter    & Leptotrichiaceae \\
\textbf{Bacteroides}       & Granulicatella  & Faecalibacterium & Lautropia        \\ 
Neisseria    & Kingella    & Campylobacter       & Actinobacillus\\
Aggregatibacter & Haemophilus   & &\\
\bottomrule
\end{tabular}
\end{table}

The results of the experiment demonstrate the efficacy of our solution sparsity in identifying key genera for clustering. Figure \ref{fig:purity} shows that the purity remains high (99\%) as $\lambda$ increases from 0 to 224 except for 121 to 210, and the number of genera used for clustering decreased from 54 to 30. Bifidobacterium, Collinsella and Escherichia leaving the former model leads to a purity drop from 99\% to 92\%. Propionibacterium leaving the former model leads to a purity drop from 99\% to 85\%. Compared to 51 core genera listed in \cite{tan2021machine} Appendix table S2, we found that using 30 core genera (Table \ref{tab:hmp}) also can reach 99\% hierarchical clustering purity. All those genera are from Bacteria kingdom. Figure \cite{fig:wob} shows the purity behavior against $\lambda$ after excluding Bacteroides. In general, purity always stays around 60\% as decreasing number of genera. The genus Bacteroides.  

\begin{figure} [!t]  
\centering
\pgfplotsset{set layers}
\begin{tikzpicture}[scale=1]
\pgfplotsset{ 
    xmin=0, xmax=300
}
\begin{axis}[ 
    xlabel={$\lambda$},
    ylabel={\ number of genera},  
    axis y line*=left, 
    ymin=0, ymax=70, 
    legend style={at={(0.45,0.6)},anchor=north east},
    ymajorgrids=true,
    grid style=dashed,
    ticklabel style = {font=\scriptsize},
] 
\addplot[line width=1pt,color=blue,mark=none]
    table [x=lambda, y=nfeatures,col sep=comma] {csv/hmp.csv}; 
\end{axis}  
\begin{axis}[  
    axis y line*=right,
    ymin=0, ymax=1,  
    yticklabel pos=right,   
    ylabel={average purity},
    ticklabel style = {font=\scriptsize},
    xticklabels={},
]
\addplot[line width=1pt,color=red,mark=none]
    table [x=lambda, y=avg_purity,col sep=comma] {csv/hmp.csv};   
\end{axis} 
\end{tikzpicture} 
\caption{The number of taxa with non-zero loadings as a function of $\lambda$ is in blue. The average purity as a function of $\lambda$ is in red. The number of taxa used for clustering is on the left y-axis for the blue line and average purity is on the right y-axis for the red line.}
\label{fig:purity}
\end{figure} 

\section{Algorithms 1 and 2 on NVIDIA Graphical Processing Units}

In this section, we discuss the implementation of Algorithms \ref{alg1} and \ref{alg2} on NVIDIA CUDA, a general-purpose parallel computing platform. We recognize that Algorithms \ref{alg1} and \ref{alg2} can be implemented in a parallel framework such as CUDA by sorting the $m$ lists independently.   

\subsection{Introduction}
Recent distributed parallel computing technologies offer a solution for handling big data by increasing overall throughput (number of jobs or tasks executed per unit of time).  The CUDA compute platform provides a scalable programming paradigm that extends C, C++, Python, and Fortran to be capable of executing parallel algorithms within thread groups on GPUs. The CUDA application is heterogeneous in the sense that parallel and sequential abstractions coexist in one application; namely, kernel and host. Users initialize a C/C++ application (hosted on a CPU) connecting to a kernel interface, which in turn allocates the resources on the GPU. Devices are responsible for performing computations and partitioning the cache. A kernel is a C++/C function with the qualifier \texttt{\_\_global\_\_}. It will run independently on GPU threads with a unique ID, as illustrated in Figure \ref{fig:cuda2}. In this toy application, 2 blocks with 4 threads are launched for parallel computation on each element of the vector, which will be cached on the GPU. The memory location of each element can be indexed with three built-in variables \texttt{threadIdx.x}, \texttt{blockIdx.x} and \texttt{blockDim.x}. We will dive into more details in the next section.
 
\subsection{Computational Speedup Results}
In this section, we run our CUDA application on a NVIDIA GeForce RTX 3060 laptop GPU with 3840 cores and 6 gigabytes of graphics memory. We run CPU implementations on an Intel 8-core I9 processor along with 40 gigabytes of memory. Figure \ref{fig:cudakernel} shows a snippet of this application that computes the quotients between the $k^{th}$ column and all other columns and stores the result in the vector $\textit{d\_out}$ through three built-in variables \texttt{threadIdx.x}, \texttt{blockIdx.x}, \texttt{blockDim.x} and \texttt{gridDim.x}. \texttt{threadIdx.x} is the index of each element in one block, and \texttt{blockIdx.x} is the index of each block in CPU memory. \texttt{gridDim.x} (the number of blocks) and \texttt{blockDim.x} (the number of threads per block) for this practical task are specified in $\lll128,128\ggg$. This tells runtime to create 128 copies of the kernel and run them in parallel. Each of these parallel invocations is a block. The code for Algorithms \ref{alg1} and \ref{alg2} is in the Appendix.

We then run this CUDA implementation with 10 replications for each size and calculate the average runtime. Table \ref{tab:speedup1} gives the speedup overview for 121 different input sizes. It shows up to 16.57x speedup over the R implementation and implies an increasing speedup as the size increases.

\begin{table}[!htb]
\centering
\caption{Speedup Results for a Matrix of Dimension Row Index $\times$ Column Header. A Value Greater than 1 Demonstrates the Efficacy of the Implementation of Algorithm \ref{alg1}.}
\csvreader[
tabular = r|rrrrrrrrrrr,
table head = &100&200&300&400&500&600&700&800&900&1000&2000 \label{tab:speedup1}\\
\hline,
late after line = \\]%
{csv/speedup.csv}{}{%
\csvcoli&\csvcolii&\csvcoliii&\csvcoliv&\csvcolv&\csvcolvi&\csvcolvii&\csvcolviii&\csvcolix&\csvcolx&\csvcolxi&\csvcolxii
}
\end{table}

\subsection{Solution Path with Varying Dimensions and $\lambda$} 
In this section, we first evaluate the behavior of the solution of Algorithm \ref{alg1} under different regularization parameters $\lambda$ and input dimensions in terms of norm $\ell^0$  and discordance. We also evaluated the space requirements for the number of breakpoints generated by Algorithm \ref{alg2}. All experiments were carried out on a CUDA GPU. 

Table \ref{tab:la1} shows the average elapsed time to compute Algorithm \ref{alg1} on 10 replications for each size. For example, Algorithm \ref{alg1} of a $5000\times1000$  matrix takes about 26 seconds and 127 seconds for a $1000\times5000$ matrix on the GPU. Since the running time of Algorithm \ref{alg1} is directly proportional to $m^2$ and $n$, the matrix of larger columns requires more computation time than the matrix of fewer columns. This can also be illustrated by the time in terms of input rows, with 5000 columns being the steepest line in Figure \ref{fig:bp2}. Lastly, Figures \ref{fig:bp1} and Table \ref{tab:breakpoints2} illustrate that the number of breakpoints generated by Algorithm \ref{alg2} is directly proportional to $n$ and $m^2$. 

\begin{table} 
\caption{Average and standard deviation time in seconds for 10 replications for each dataset with varying number of columns with fixed number of rows at 1000, 2000, and 5000. }
\label{tab:la1}  
\centering
\csvreader[
tabular = r|rrr,
table head =  & 1000 & 2000 & 5000\\
\hline,
late after line = \\]%
{csv/sparsel1.csv}{}{%
\csvcoli&\csvcolv~$_{\csvcolxxiii}$&\csvcolvi~$_{\csvcolxxiv}$ & \csvcolvii~$_{\csvcolxxv}$
}  
\end{table}
  
\begin{table} 
\caption{Average and standard deviation time in seconds over 10 replications for each dataset with varying the number of rows with the fixed number of columns at 1000,2000, and 5000. }
\label{tab:la2}   
\centering
\csvreader[
tabular = r|rrr,
table head =  & 1000 & 2000 & 5000\\
\hline,
late after line = \\]%
{csv/sparsel1.csv}{}{%
\csvcoli&\csvcolii~$_{\csvcolxx}$&\csvcoliii~$_{\csvcolxxi}$ & \csvcoliv~$_{\csvcolxxii}$
} 
\end{table}

\begin{figure*}[!t]
\centering
\subfloat[]{
\begin{tikzpicture}[scale=1]
\begin{axis}[ 
    xlabel={Input columns},  
    xmin=100, xmax=1000,
    ylabel={Running time(seconds)},
    ymin=0, ymax=130, 
    legend pos= north west,
    ymajorgrids=true,
    grid style=dashed,
    axis line style={gray},
    legend style={draw=none},
    ticklabel style = {font=\scriptsize},
]
\addplot[
    color=blue,
    mark=square,
    ]
    table [x=x, y=n1000,col sep=comma] {csv/sparsel1.csv};
    \legend{n=1000}
\addplot[
    color=red,
    mark=*,
    ]
    table [x=x, y=n2000,col sep=comma] {csv/sparsel1.csv};
    \addlegendentry{n=2000}
\addplot[
    color=green,
    mark=triangle,
    ]
    table [x=x, y=n5000,col sep=comma] {csv/sparsel1.csv};
    \addlegendentry{n=5000}
\end{axis}
\end{tikzpicture}}
\subfloat[]{
\begin{tikzpicture}[scale=1]
\begin{axis}[
    xlabel={Input rows},
    xmin=100, xmax=1000,
    axis line style={gray},
    legend style={draw=none},
    ymin=0, ymax=130, 
    legend pos= north west,
    ymajorgrids=true,
    grid style=dashed,
    ticklabel style = {font=\scriptsize},
]
\addplot[
    color=blue,
    mark=square,
    ]
    table [x=x, y=m1000,col sep=comma] {csv/sparsel1.csv};
    \legend{m=1000}
\addplot[
    color=red,
    mark=*,
    ]
    table [x=x, y=m2000,col sep=comma] {csv/sparsel1.csv};
    \addlegendentry{m=2000}
\addplot[
    color=green,
    mark=triangle,
    ]
    table [x=x, y=m5000,col sep=comma] {csv/sparsel1.csv};
    \addlegendentry{m=5000}
\end{axis}
\end{tikzpicture}}
\caption{Running time for Algorithm \ref{alg2}.}
\label{fig:bp2}
\end{figure*}
\begin{figure*} [!t]
\centering 
\subfloat[]{\begin{tikzpicture}[scale=1]
\begin{axis}[
    xlabel={Input rows},
    ylabel={Breakpoints(millions)},
    xmin=100, xmax=1000,
    ymin=0, ymax=1300, 
    legend pos= north west,
    ymajorgrids=true,
    grid style=dashed,
    axis line style={gray},
    legend style={draw=none},
    ticklabel style = {font=\scriptsize},
]
\addplot[
    color=blue,
    mark=square,
    ]
    table [x=x, y=lm1000,col sep=comma] {csv/sparsel1.csv};
    \legend{m=1000}
\addplot[
    color=red,
    mark=*,
    ]
    table [x=x, y=lm2000,col sep=comma] {csv/sparsel1.csv};
    \addlegendentry{m=2000}
\addplot[
    color=green,
    mark=triangle,
    ]
    table [x=x, y=lm5000,col sep=comma] {csv/sparsel1.csv};
    \addlegendentry{m=5000}
\end{axis}
\end{tikzpicture}}
\subfloat[]{
\begin{tikzpicture}[scale=1]
\begin{axis}[ 
    xlabel={Input columns}, 
    xmin=100, xmax=1000,
    ymin=0, ymax=1300, 
    legend pos= north west,
    ymajorgrids=true,
    grid style=dashed,
    axis line style={gray},
    legend style={draw=none},
    ticklabel style = {font=\scriptsize},
]
\addplot[
    color=blue,
    mark=square,
    ]
    table [x=x, y=ln1000,col sep=comma] {csv/sparsel1.csv};
    \legend{n=1000}
\addplot[
    color=red,
    mark=*,
    ]
    table [x=x, y=ln2000,col sep=comma] {csv/sparsel1.csv};
    \addlegendentry{n=2000}
\addplot[
    color=green,
    mark=triangle,
    ]
    table [x=x, y=ln5000,col sep=comma] {csv/sparsel1.csv};
    \addlegendentry{n=5000}
\end{axis}
\end{tikzpicture}}
\caption{Number of breakpoints for Algorithm \ref{alg2}.}
\label{fig:bp1}
\end{figure*}

\begin{table} 
\caption{Average and standard deviation of breakpoints in millions over 10 replications varying the number of columns with the fixed number of rows at 1000, 2000, and 5000..}
\centering 
\csvreader[
tabular = r|rrr,
table head =  & 1000 & 2000 & 5000 \label{tab:breakpoints2}\\
\hline,
late after line = \\]%
{csv/sparsel1.csv}{}{%
\csvcoli&\csvcolxvii~$_{\csvcolxxxv}$&\csvcolxviii~$_{\csvcolxxxvi}$ & \csvcolxix~$_{\csvcolxxxvii}$
} 
\end{table}

\begin{table} 
\caption{Average and standard deviation of breakpoints in millions over 10 replications varying the number of rows with the fixed number of columns at 1000, 2000, and 5000.}
\centering 
\csvreader[
tabular = r|rrr,
table head =  & 1000 & 2000 & 5000\\
\hline,
late after line = \\]%
{csv/sparsel1.csv}{}{%
\csvcoli&\csvcolxiv~$_{\csvcolxxxii}$&\csvcolxv~$_{\csvcolxxxiii}$ & \csvcolxvi~$_{\csvcolxxxiv}$
}
\end{table}
  
\section{Conclusion}
We have introduced a new and efficient method to estimate a sparse and outlier-resistant best-fit line.  The development is based in part on linear programming theory. We demonstrate that the relevant LP can be solved via sorting lists of numbers.  In addition, our algorithms can be processed in parallel for increased efficiency, which may enable a wide range of new practical applications.  The new method is compared to state-of-the-art methods on synthetic datasets and is shown to provide insights with human micorbiome data.

\section*{Acknowledgments}
This should be a simple paragraph before the references to thank the individuals and institutions who have supported your work on this article.

\newpage
 
\bibliographystyle{unsrt}  
\bibliography{ref}

\end{document}